\documentclass[12pt]{amsart}
\usepackage{graphicx}
\usepackage{amssymb}
\usepackage{amsfonts}
\usepackage{amsmath}
\usepackage{array}
\usepackage{rotating}

\usepackage{pstricks}
\usepackage{array,delarray}

\usepackage{placeins}

\usepackage{pspicture}
\usepackage[matrix,ps,xdvi]{xypic} % ne pas enlever
\newdimen\vcadre\vcadre=0.1cm % marges verticales de la boite
\newdimen\hcadre\hcadre=0.1cm % marges horizontales de la boite
\def\GrTeXBox#1{\vbox{\vskip\vcadre\hbox{\hskip\hcadre%
      $#1$%
   \hskip\hcadre}\vskip\vcadre}}
\def\arx#1[#2]{\ifcase#1 \relax \or%
  \ar @{-}[#2]  \or%
  \ar @2{-}[#2] \or%
  \ar @{--}[#2] \or%
  \ar @2{.}[#2] \or%
  \ar @{~}[#2]  \fi}

\def\arbg#1#2{
\newdimen\vcadre\vcadre=0.01cm % marges verticales de la boite
\newdimen\hcadre\hcadre=0.01cm % marges horizontales de la boite
\xymatrix@R=0.1cm@C=2mm{
 && {\GrTeXBox{#1}}\arx1[ld]\arx1[rd]\\
& {\GrTeXBox{#2}}\arx1[dl]\arx1[dr] && {\GrTeXBox{z}}\\
 {\GrTeXBox{x}} &&  {\GrTeXBox{y}}\\
}
}
\def\arbd#1#2{
\newdimen\vcadre\vcadre=0.01cm % marges verticales de la boite
\newdimen\hcadre\hcadre=0.01cm % marges horizontales de la boite
\xymatrix@R=0.1cm@C=2mm{
 & {\GrTeXBox{#1}}\arx1[ld]\arx1[rd]\\
{\GrTeXBox{x}} & & {\GrTeXBox{#2}}\arx1[ld]\arx1[rd]\\
& {\GrTeXBox{y}}&  & {\GrTeXBox{z}}\\
}
}
\def\arbgnew#1#2{
\newdimen\vcadre\vcadre=0.01cm % marges verticales de la boite
\newdimen\hcadre\hcadre=0.01cm % marges horizontales de la boite
\xymatrix@R=0.1cm@C=2mm{
 && {\GrTeXBox{#1}}\arx1[ld]\arx1[rd]\\
& {\GrTeXBox{#2}}\arx1[dl]\arx1[dr] && {\GrTeXBox{B}}\\
 {\GrTeXBox{.}} &&  {\GrTeXBox{A}}\\
}
}

\newtheorem{example}{Example}[section]

\newtheorem{theorem}[example]{Theorem}

\newtheorem{corollary}[example]{Corollary}

\newtheorem{definition}[example]{Definition}
\newtheorem{proposition}[example]{Proposition}
\newtheorem{algorithm}[example]{Algorithm}
\newtheorem{lemma}[example]{Lemma}

\def\Proof{\noindent \it Proof -- \rm}
\def\qed{\hspace{3.5mm} \hfill \vbox{\hrule height 3pt depth 2 pt width 2mm}
\bigskip}

\def\gaudend{\prec}     % gauche dendriforme
\def\droitdend{\succ}   % droite dendriforme
\def\gautrid{\!\prec\!}   % gauche trigebre dendriforme
\def\miltrid{\circ}       % milieu trigebre dendriforme
\def\droittrid{\!\succ\!} % droite trigebre dendriforme

\def\TD{{\mathfrak TD}}
\def\convW{{*_W}}
\def\K{{\mathbb K}}
\def\PW{{\rm PW}}
\def\pack{{\rm pack}}
\def\PBT{{\bf PBT}}

\def\maj{{\rm maj}}
\def\FQSym{{\bf FQSym}}

\def\PQSym{{\bf PQSym}}
\def\CQSym{{\bf CQSym}}
\def\SQSym{{\bf SQSym}}
\def\WQSym{{\bf WQSym}}
\def\ev{{\rm ev}}
\def\pev{{\rm pev}}
\def\qrpark{{\bf q}}
\def\q{{\bf q}}

\def\ssh{\Cup}

\def\sconc{\bullet}
\def\std{{\rm Std}}
\def\std{{\rm std}}
\def\Park{{\rm park}}

\def\park{{\bf a}}
\def\a{{\bf a}}
\def\b{{\bf b}}
\def\c{{\bf c}}

\def\<{\langle}
\def\>{\rangle}

\def\Z{\operatorname{\mathbb Z}}
\def\N{\operatorname{\mathbb N}}
\def\KK{\operatorname{\mathbb K}}

\def\F{{\bf F}}

\def\T{{\bf T}}
\def\G{{\bf G}}

\def\M{{\bf M}}
\def\P{{\bf P}}

\def\SG{{\mathfrak S}}

\def\A{{\sf A}}

\def\Sym{{\bf Sym}}

\def\NDPF{{\rm NDPF}}

\def\Q{{\bf Q}}

\def\sinv{{\rm sinv}}
\def\smaj{{\rm smaj}}

\def\PF{{\rm PF}}

\def\D{{\mathcal D}}
\def\sep{|}
\def\qrpark{{\bf q}} 
\def\park{{\bf a}}    % une fonction de parking
\def\PS{{\mathcal P}} 

\def\shuff#1#2{\mathbin{
\hbox{\vbox{ \hbox{\vrule \hskip#2 \vrule height#1 width 0pt
}%
\hrule}%
\vbox{ \hbox{\vrule \hskip#2 \vrule height#1 width 0pt
\vrule }%
\hrule}%
}}}

\def\shuf{{\mathchoice{\shuff{7pt}{3.5pt}}%
{\shuff{6pt}{3pt}}%
{\shuff{4pt}{2pt}}%
{\shuff{3pt}{1.5pt}}}}%
\def\shuffle{\,\shuf\,}

\def\Tabvrule{\vrule width-0.4pt}       % Difference de largeur
\def\Tabhrule{\hrule \hrule height-0.4pt} % Difference de hauteur
\def\Tabstrut{\vrule height2.2ex % Sur la ligne
                     depth0.8ex  % Sous la ligne
                     width0ex    % centrage horizontal
\relax}

\def\PasCase#1{\omit%
            $\vcenter{\hbox {\vbox to 0.4pt{}}
               \hbox{\makebox[3ex]{\Tabstrut$#1$}}}%
               \Tabvrule$}
\def\PasCasePoint{\PasCase{\cdot}}
\def\DessinCarre#1{%
    \vcenter{\hbox{}\hrule
             \hbox{\vrule\makebox[3ex]{\Tabstrut$#1$}\vrule}\Tabhrule}%
             \Tabvrule}
\def\GenRuban#1{\vcenter{\halign{&$\DessinCarre{##}$\cr#1}}\egroup}

\def\sTabvrule{\vrule width-0.4pt}
\def\sTabhrule{\hrule \hrule height-0.4pt}
\def\sTabstrut{\vrule height1.6ex depth0.6ex width0ex \relax}
\def\sDessinCarre#1{%
    \vcenter{\hbox{}\hrule
             \hbox{\vrule\makebox[2.3ex]%
                  {\sTabstrut$\scriptstyle#1$}\vrule}\sTabhrule}%
             \sTabvrule}
\def\sGenRuban#1{\vcenter{\halign{&$\sDessinCarre{##}$\cr#1}}\egroup}

\def\ruban{%
  \bgroup
  \let\ =\omit
  \let\\=\cr
  \let\x=\times
  \let\.=\PasCasePoint
  \offinterlineskip
  \GenRuban}

\def\sruban{%
  \bgroup
  \let\ =\omit
  \let\x=\times
  \let\\=\cr
  \offinterlineskip
  \sGenRuban}

\title{Duplicial algebras, parking functions, and Lagrange inversion}
\author[J.-C.~Novelli and J.-Y.~Thibon]%
{Jean-Christophe Novelli and Jean-Yves Thibon}

\address[] {Institut Gaspard Monge, Universit\'e de Marne-la-Vall\'ee \\
5 Boulevard Descartes \\Champs-sur-Marne \\77454 Marne-la-Vall\'ee cedex 2 \\
FRANCE}
\email[Jean-Christophe Novelli]{novelli@univ-mlv.fr}
\email[Jean-Yves Thibon]{jyt@univ-mlv.fr} 
\date{}

\keywords{Operads, Noncommutative symmetric functions,
parking functions, Lagrange inversion
}
\subjclass{18D50,05E05,16T30}

\begin{document}

\begin{abstract}
We provide operadic interpretations for two Hopf subalgebras of the algebra of
parking functions. The Catalan subalgebra is identified with the free
duplicial algebra on one generator, and the Schr\"oder subalgebra is
interpreted by means of a new operad, which we call triduplicial.

The noncommutative Lagrange inversion formula is then interpreted in terms of
duplicial structures.
The generic solution of the noncommutative inversion problem appears  as the
formal sum of all parking functions. This suggests that combinatorial
generating functions derived by functional inversion should be obtainable by
evaluating a suitable character on this generic solution.  This idea is
illustrated by means of the Narayana polynomials, of which we obtain bivariate
``super-analogues'' by lifting to parking functions a classical character of
the algebra of symmetric functions.  Other characters, such as evaluation of
symmetric functions on a binomial element, are also discussed.
\end{abstract}

\maketitle
\footnotesize
\tableofcontents
\normalsize
%%%%%%%%%%%%%%%%%%%%%%%%%%%%%%%%%%%%%%%%%%%%%%%%%%%%%%%%%%%%%%%%%%%%%%%%%%%%%%%
%%%%%%%%%%%%%%%%%%%%%%%%%%%%%%%%%%%%%%%%%%%%%%%%%%%%%%%%%%%%%%%%%%%%%%%%%%%%%%%
%%%%%%%%%%%%%%%%%%%%%%%%%%%%%%%%%%%%%%%%%%%%%%%%%%%%%%%%%%%%%%%%%%%%%%%%%%%%%%%
\section{Introduction}

In its simplest version, the Lagrange inversion formula gives the coefficients
of the unique formal power series 
\begin{equation}
g(t)=\sum_{n\geq0}g_nt^{n+1}
\end{equation} 
satisfying the functional equation
\begin{equation}\label{eqlag}
t = \frac{g}{\varphi(g)} \,
\end{equation}
where
\begin{equation}
\varphi(x)=\sum_{n\ge 0}a_nx^n\quad (a_0\not=0)\,
\end{equation}
is a series with nonzero constant term.

As is well-known ({\it cf.} \cite[5.4]{Stan2}), the result is 
\begin{equation}\label{lag1}
g_n =\frac1{n+1}[x^n](\varphi(x))^{n+1}\,.
\end{equation}

There exist several combinatorial variations of this result,
providing descriptions of $g_n$ as a polynomial in the $a_n$
as well as $q$-analogues or noncommutative versions.

As any generic identity among formal power series, the formula of Lagrange is
equivalent to an identity on symmetric functions \cite{Las}. One may identify the $a_n$
with a sequence of algebraically independent generators of the algebra
$Sym(X)$ of symmetric functions. For example, set $a_n=h_n(X)$
(the complete homogeneous functions as
in~\cite[Ex. 24 p. 35, Ex. 25 p. 132]{Mcd})
so that
\begin{equation}
\varphi(u)=\sum_{n\ge 0}h_n(X)u^n= \prod_{n\ge 1}(1-ux_n)^{-1}=: \sigma_u(X)
\end{equation}
and the result \eqref{lag1} reads (in $\lambda$-ring notation)
\begin{equation}\label{eq:comLag}
g_n =\frac1{n+1}h_n((n+1)X)\,.
\end{equation}
What makes this expression interesting is that this symmetric
function is the Frobenius characteristic of the permutation
representation of the symmetric group $\SG_n$ on the set $\PF_n$ of
parking functions of length $n$  \cite{Hai1}.

It has moreover been shown in \cite{NTLag} that  noncommutative
versions of the Lagrange formula such as in \cite{Ges} or \cite{PPR}
could be formulated in terms of noncommutative symmetric functions~\cite{NCSF1}.
Again, the term $g_n$ of degree $n$ in the series $g$ satisfying
\begin{equation}
g = \sum_{n\ge 0} S_n g^n
\end{equation}
(where the $S_n$ are the noncommutative complete symmetric functions)
is the Frobenius characteristic of the natural representation of the $0$-Hecke
algebra $H_n(0)$ on parking functions \cite{NTLag}. 

Now, there is a combinatorial Hopf algebra, $\PQSym$, which is based
on parking functions. Denoting its natural basis by $\F_\a$ as in \cite{NT1},
we shall see that the series $g$ is the image of the sum of all parking
functions
\begin{equation}
G = \sum_{\a\in\PF} \F_\a
\end{equation}
by a homomorphism of Hopf algebras. This morphism can be defined at the
level of several intermediate Hopf algebras, such as  $\SQSym$
or $\CQSym$ (both defined in \cite{NTPark}), 
which have as graded dimensions the sequences of little Schr\"oder
numbers and Catalan numbers respectively. We shall see that $\CQSym$
can be naturally identified with the free duplicial algebra on one
generator, and that as an element of $\CQSym$, $G$ satisfies a quadratic
functional equation allowing to identify it with the formal sum
of all binary trees. This forces a bijection between these trees and
nondecreasing parking functions, which in turn implies the existence of
an involution explaining the symmetry properties of $G$ observed
in \cite{NTLag} and established there by means of a different involution.

The operadic interpretation of $\CQSym$ suggests the existence of a similar
one for $\SQSym$. We thus introduce a new operad, which we
call triduplicial, for which $\SQSym$ is the free algebra on one generator.

Finally, we illustrate the idea that combinatorial series derived by
means of Lagrange inversion should come from some character of the algebra
of parking functions. First, we generalize  Lassalle's expression \cite{Lass} of the
Narayana polynomials $c_n(t)$ to super-Narayana polynomials $P_n(t,q)$
counting signed parking functions according to certain statistics, and such
that $P_n(t,0)=(1+t)c_n(1+t)$. This last polynomial is known to count
Schr\"oder paths according to the number of horizontal steps
\cite[A060693]{Slo}.

These paths can be naturally encoded by a subset of signed parking functions,
such that minus signs correspond to horizontal steps. 
We conclude by investigating the character of symmetric functions defined
by evaluation on a binomial element, and obtain a combinatorial interpretation
of its lift to parking functions.

\medskip
This paper is  a continuation of~\cite{NTPark, NTLag}. We have only 
recalled the most basic definitions so as to make it reasonably self-contained.

%%%%%%%%%%%%%%%%%%%%%%%%%%%%%%%%%%%%%%%%%%%%%%%%%%%%%%%%%%%%%%%%%%%%%%%%%%%%%%%
%%%%%%%%%%%%%%%%%%%%%%%%%%%%%%%%%%%%%%%%%%%%%%%%%%%%%%%%%%%%%%%%%%%%%%%%%%%%%%%
%%%%%%%%%%%%%%%%%%%%%%%%%%%%%%%%%%%%%%%%%%%%%%%%%%%%%%%%%%%%%%%%%%%%%%%%%%%%%%%
\section{Operads and combinatorial Hopf algebras}

Combinatorial Hopf algebras are certain graded bialgebras based on
combinatorial objects. This is a heuristic concept, and there is no general
agreement on what should be their general definition. For us, they arise as
natural generalizations of the algebra of symmetric functions. Surprisingly,
most of these algebras also arise in the theory of operads (sometimes with
non-obviously equivalent definitions). We shall see that our analysis of the
noncommutative Lagrange inversion problem will allow to identify the operads
associated with two algebras ($\CQSym$ and $\SQSym$) from \cite{NTPark}, for
which the operadic interpretation was not known.

%%%%%%%%%%%%%%%%%%%%%%%%%%%%%%%%%%%%%%%%%%%%%%%%%%%%%%%%%%%%%%%%%%%%%%%%%%%%%%%
\subsection{Noncommutative symmetric functions}

The Hopf algebra of noncommutative symmetric functions~\cite{NCSF1},
denoted by $\Sym$, or by $\Sym(A)$ if we consider the realization in terms of
an auxiliary alphabet, is the free associative algebra over a sequence
$(S_n)_{n\ge 1}$, graded by ${\rm deg}\,(S_n)=n$.
For $A=\{a_n|n\ge 1\}$ a totally ordered set of noncommuting indeterminates,
we define the generating series
\begin{equation}
\sigma_t(A):=\sum_{n\ge 0}S_n(A)t^n =\prod_{n\ge 1}^\rightarrow (1-ta_i)^{-1}
\end{equation}
so that $S_n(A)$ is the sum of all nondecreasing words of length $n$
(and $S_0=1)$.

Bases of $\Sym_n$ are labelled by compositions $I=(i_1,\ldots,i_r)$ of $n$. 
The natural basis is $S^I=S_{i_1}\cdots S_{i_r}$. 
The reverse refinement order $I\le J$ on compositions of $n$ means that
the parts of $I$ are sums of consecutive parts of $J$.
The length $r$ of $I$ is denoted by $\ell(I)$.
The conjugate composition is denoted by $I^\sim$. 
The ribbon basis is defined by
\begin{equation}
S^I = \sum_{J\le I}R_J \quad \Longleftrightarrow R_I
    =\sum_{J\le I}(-1)^{\ell(I)-\ell(J)}S^J\,.
\end{equation}
For two compositions $I=(i_1,\ldots,i_r)$ and $J=(j_1,\ldots,j_s)$, define
\begin{equation}
I\cdot J = (i_1,\ldots,i_r,j_1,\ldots,j_s)\
\text{and}\
I \triangleright J = (i_1,\ldots,i_r+j_1,\ldots,j_s)\,.
\end{equation}
Then,
\begin{equation}
R_I R_J = R_{I\cdot J} + R_{I\triangleright J}\,.
\end{equation}
Thus, in the $R$-basis, the multiplication of $\Sym$ is the sum of the two
operations $\cdot$ and $\triangleright$. For this reason, one may consider
$\Sym$ as a \emph{dialgebra} \cite{Lod}. Both operations are associative, and
furthermore
\begin{equation}
(I *_1 J) *_2 K = I *_1 (J *_2 K)
\end{equation}
where $*_1$ and $*_2$ are either $\cdot$ or $\triangleright$.
Such an algebra is called an $As^{(2)}$-algebra \cite{Zin}.
It is easy to see that $\Sym$ is actually the free  $As^{(2)}$-algebra
on one generator. This is the simplest example of
a combinatorial Hopf algebra associated with an operad.

A somewhat more interesting example arises if we consider a slightly more
general version of the  noncommutative Lagrange inversion problem.

%%%%%%%%%%%%%%%%%%%%%%%%%%%%%%%%%%%%%%%%%%%%%%%%%%%%%%%%%%%%%%%%%%%%%%%%%%%%%%%
\subsection{Nondecreasing parking functions}

The versions of \cite{Ges} and \cite{PPR} of the noncommutative
inversion formula can be interpreted as solving the equation
\begin{equation}
\label{genfeq}
g = S_0 +S_1 g + S_2g^2+S_3g^3+\cdots\,,
\end{equation}
where $S_0$ is now another indeterminate not commuting
with the other ones \cite{NTLag}. 

A \emph{nondecreasing parking function} is a 
nondecreasing word $\pi=\pi_1\ldots\pi_n$
over $[n]$ such that $\pi_i\le i$. A \emph{parking
function} is any rearrangement of such a word.

The solution of \eqref{genfeq} can be expressed
\cite{NTLag}
 in the form\footnote{
Although Equation \eqref{eq:comLag} would make sense for noncommutative symmetric
functions, its does not hold in this context: one can compute
\begin{equation}
\frac14 S_3(4A)=S_3+\frac32 S^{21}+\frac32 S^{12}+S^{111}
\not=g_3=S_3+2S^{21}+S^{12}+S^{111}.
\end{equation}
However, it is possible to derive \eqref{eq:comLag} from \eqref{gensol}.
This is essentially the classical combinatorial proof, originally due to Raney \cite{Ran}.
Here is a somewhat simplified version of the argument.

For a word $w$, let $\lambda(w)$ be the partition obtained by removing the zeros
in the evaluation of $w$ and sorting the parts in decreasing order.
Projecting \eqref{gensol} by $S_k\mapsto h_k(X)$, we see that
the coefficient of $h_\mu$ in the $g_n$ of \eqref{eq:comLag}
 is the number of nondecreasing parking functions
$\pi\in\NDPF_n$ such that $\lambda(\pi)=\mu$.

This number can be obtained as follows. The words encoding the terms in the $g_n$ of \eqref{gensol}
are the evaluations of $\NDPF_n$ with an extra 0, hence of sum $n$ and length $(n+1)$.
The sum of the $S^w$ for of all such words $w$ is the term of degree $n$ in
\begin{equation}
(S_0+S_1+S_2+\cdots +S_n+\cdots)^{n+1}
\end{equation}
and under the substitution $S_k = h_k(X)$, this sum becomes $h_n((n+1)X)$.

Now, if $w$ is a word on $\N$ of sum $n$ and length $(n+1)$, there is
exactly one cyclic shift of $w$ which is the evaluation of a nondecreasing
parking function followed by a 0, or, equivalently, the Polish code of a rooted plane tree with $n+1$ vertices
(this is the so-called {\it cycle lemma}).

For example, among the seven words
\begin{equation}
3010200,\ 0102003,\ 1020030,\ 0200301,\ 2003010,\ 0030102,\ 0301020
\end{equation}
only the first one is a Polish code.

Thus, $h_n((n+1)X)=(n+1)g_n$.

Here is a simple proof of the cycle lemma. A word of sum $n$ and of length $n+1$ cannot be a nontrivial
power of another word, hence has $n+1$ distinct cyclic shifts.
The number of such words is the coefficient of $x^n$ in
$$
(x^0+x^1+x^2+\cdots)^{n+1}=(1-x)^{-n-1}=\sum_{m\ge 0}{m+n-1\choose m}x^m,
$$
that is, ${2n\choose n}$. Thus, the number of cyclic orbits is the Catalan
number $\frac1{n+1}{2n\choose n}$, which is precisely the number of rooted plane trees with $n+1$ vertices.
Hence, each orbit contains exactly one tree.
}
\begin{equation}\label{gensol}
g_n = \sum_{\pi\in\NDPF_n} S^{\ev(\pi)\cdot 0}\,,
\end{equation}
where $\NDPF_n$ denotes the set of nondecreasing parking functions
of length $n$. The \emph{evaluation} of a word is the sequence 
$\ev(w)=(|w|_i)$ recording the number of occurences of each letter $i$.
For example,
\begin{eqnarray}
\label{sol1110}
g_0=S_0\,,\ g_1= S_1S_0=S^{10}\,,\ g_2=S^{110}+S^{200}\,,\\
g_3=S^{1110}+S^{1200}+S^{2010}+S^{2100}+S^{3000}\,,
\end{eqnarray}
the nondecreasing parking functions giving $g_3$ being
(in this order) $123$, $122$, $113$, $112$, $111$.
%%%

%%%

As shown in \cite{NTPark}, there is a combinatorial Hopf algebra
based on nondecreasing parking functions. We shall see  that
it is directly related to Lagrange inversion, and uncover
its operadic interpretation.

%%%%%%%%%%%%%%%%%%%%%%%%%%%%%%%%%%%%%%%%%%%%%%%%%%%%%%%%%%%%%%%%%%%%%%%%%%%%%%%
\subsection{Duplicial algebras}

A duplicial algebra \cite{LodTrip} is a vector space endowed with two
bilinear associative operations $\prec$ and $\succ$ such that
\begin{equation}
(x\succ y)\prec z = x\succ(y\prec z)\,.
\end{equation}
It is known that the free duplicial algebra $\D$  on one generator has
a basis labelled by binary trees. The operations  $\prec$ and $\succ$ 
can be respectively identified with the products $\backslash$ (under) and
$/$ (over) (see \cite{LodTrip}). The dimensions of the homogeneous components
are therefore the Catalan numbers $1,2,5,14\dots$. 

Here is another realization of $\D$.

The Hopf algebra $\PQSym$, defined in \cite{NT1}, is the linear span
of elements $\F_\a$ where $\a$ runs over all parking functions.
%%%
Its product rule is given by the shifted shuffle\footnote{
For a word $w$ on the alphabet $\{1,2,\ldots\}$, denote by $w[k]$ the word
obtained by replacing each letter $i$ by the integer $i+k$.
If $u$ and $v$ are two words, with $u$ of length $k$, one defines
the {\em shifted concatenation}
\begin{equation}
u\sconc v = u\cdot (v[k])
\end{equation}
and the {\em shifted shuffle}
\begin{equation}
u\ssh v= u\shuffle (v[k])\,.
\end{equation}
where $\shuffle$ is the usual shuffle product on words defined recursively
by
\begin{equation}
(au)\shuffle (bv) = a\cdot(u\shuffle (bv)) + b\cdot ((au)\shuffle v),
\end{equation}
with $u\shuffle\epsilon=\epsilon\shuffle u=u$ if $\epsilon$ is the empty word.
}
\begin{equation}
\label{prodF}
\F_{\park'}\F_{\park''}:=\sum_{\park\in\park'\ssh\park''}\F_\park\,.
\end{equation}
For example,
\begin{equation}
\F_{12}\F_{11}= \F_{1233} + \F_{1323} + \F_{1332} + \F_{3123} + \F_{3132}
+ \F_{3312}\,.
\end{equation}
The coproduct of $\PQSym$ is defined in terms of the operation of
{\it parkization}, defined in \cite{NTPark}.
This operation associates with a word  $w$ on $\{1,2,\ldots\}$ a
parking function $\Park(w)$, which coincides with the standardization
$\std(w)$ when $w$ is a word without repeated letters. It is computed
by the following 

{\footnotesize
\begin{algorithm}
\label{parkisation}
\noindent
\emph{Input}:  A word $w$.

\noindent
\emph{Output}:  A parking function $\Park(w)$.

Let $n$ be the length of $w$.
Define
\begin{equation}
\label{dw}
d(w):=\min \{i\ |\ \sharp|\{w_j\leq i\}<i \}\,.
\end{equation}

\begin{itemize}
\item If $d(w)=n+1$, return $w$.
\item Otherwise, let $w'$ be the word obtained by decrementing all
the elements of $w$ greater than $d(w)$. Then return the parkized word
of $w'$.
\end{itemize}
\end{algorithm}
}

The coproduct
\begin{equation}
\label{coprod}
\Delta \F_{\park}:= \sum_{u\cdot v=\park} \F_{\Park(u)} \otimes \F_{\Park(v)}.
\end{equation}
endows $\PQSym$ with the structure of a graded connected bialgebra, hence of a Hopf algebra.

For example,
\begin{equation}
\Delta\F_{3132} = 1\otimes\F_{3132} + \F_{1}\otimes\F_{132} +
\F_{21}\otimes\F_{21} + \F_{212}\otimes\F_{1} + \F_{3132}\otimes 1\,.
\end{equation}

%%%

\bigskip
The sums
\begin{equation}
\P^\pi := \sum_{\park ; \park^\uparrow=\pi} {\F_\park}
\end{equation}
where $\park^\uparrow$ means the non-decreasing reordering and $\pi$ runs over
the set $\NDPF$ of non-decreasing parking functions, span a cocommutative Hopf
subalgebra $\CQSym$ of $\PQSym$ \cite{NTPark}.

The basis $\P^\pi$ is multiplicative:
\begin{equation}
\P^\alpha\P^\beta=\P^{\alpha\bullet \beta}
\end{equation}
where $\alpha\bullet\beta$ is the usual shifted concatenation
$\alpha\cdot\beta[k]$ if $\alpha$ is of length $k$. For example,
$\P^{12}\P^{113}=\P^{12335}$.

The $S_n=\P^{1^n}$ generate a Hopf subalgebra isomorphic to $\Sym$,
which is also a quotient of $\CQSym$ (see below).

Introduce now a second product involving a different kind of shifted
concatenation:
\begin{equation}
\P^\alpha\prec\P^\beta
=\P^{\alpha\cdot\beta[\max(\alpha)-1]}=:\P^{\alpha\circ\beta}\,.
\end{equation}
For example, $\P^{12}\prec\P^{113}=\P^{12224}$.

\begin{proposition}
\label{prop-cqs}
Let us denote the usual product of $\CQSym$ by $\succ$.

\noindent
Then $(\CQSym,\prec,\succ)$ is the free (unitary) duplicial algebra on one generator
$x=\P^1$.
\end{proposition}

\Proof It is clear that $\succ$ and $\prec$ are associative. 
Also, 
\begin{equation}
(\P^{\alpha}\succ \P^\beta)\prec \P^\gamma =
\P^{\alpha}\succ (\P^\beta\prec \P^\gamma)
\end{equation}
since, denoting by $n$ the length of $\alpha$ and
by $m$ the maximum letter of $\beta$,
\begin{equation}
(\alpha\bullet\beta)\circ\gamma= \alpha\cdot\beta[n]\cdot\gamma[m+n-1]
=\alpha\bullet(\beta\circ\gamma)\,.
\end{equation}

Thus,  $(\CQSym,\prec,\succ)$ is duplicial.

The other cross-associativity relation is not satisfied, since already
$(1\circ 1)\bullet 1=113\not= 112=1\circ (1\bullet 1)$.

Actually, $\CQSym$ is free on the generator
$x=\P^1$, since all non decreasing parking functions can be obtained by
iterating $\bullet$ and $\circ$ on $1$. Indeed, $\pi\in\NDPF$ is either
non-connected, that is, of the form $\pi=\pi'\bullet \pi''$ with nontrivial
factor, or connected,
which means that $\pi_i<i$ for all $i>1$,
so that it can be written $\pi=1\pi'=1\circ\pi'$
for some other $\pi'$. Since nondecreasing parking functions
of length $n$ are in bijection with
binary trees on $n$ nodes, we see that $\CQSym$ is indeed free. \qed 

We shall see in Section \ref{sec:invol} that $\CQSym$ has an involution $\iota$
exchanging the two products. 

%%%%%%%%%%%%%%%%%%%%%%%%%%%%%%%%%%%%%%%%%%%%%%%%%%%%%%%%%%%%%%%%%%%%%%%%%%%%%%%
\subsection{The duplicial operad}

The duplicial operad is described in detail in \cite{LodTrip}. Here is a brief
summary.
First, there is a notion of duplicial coproduct, allowing to define
primitive elements.

A duplicial bialgebra \cite{LodTrip} is a duplicial algebra endowed with
a coproduct $\delta$ satifsfying 
\begin{equation}\label{codup}
\delta(x*y)
= x\otimes y + \sum_{(x)}x_{(1)}\otimes(x_{(2)}*y) +
  \sum_{(y)}(x*y_{(1)})\otimes y_{(2)},
\end{equation}
where
$\delta x=\sum_{(x)}x_{(1)}\otimes x_{(2)}$ (Sweedler's notation) and
$*$ is $\prec$ or $\succ$.

The usual coproduct of $\CQSym$ is defined in terms of the parkization
operation 
%%%
%
\begin{equation}
\label{coprodP}
\Delta\P^{\pi} = \sum_{u,v ; (u.v)^\uparrow=\pi}
{\P^{\Park(u)} \otimes \P^{\Park(v)}}\,,
\end{equation}
where $u$ and $v$ run over the set of non-decreasing words.
%%%
%
Selecting the terms coming from the deconcatenations of $\pi$ in this coproduct, we obtain a duplicial coproduct: 
\begin{proposition}
The coproduct
\begin{equation}
\delta(\P^\pi)
=\sum_{k=1}^{n-1}\P^{\Park(\pi_1 \cdots\pi_k)}\otimes \P^{\Park(\pi_{k+1}\cdots\pi_n)} 
\end{equation}
endows $\CQSym$ with the structure of a duplicial bialgebra. 
\end{proposition}

\Proof
An immediate verification.
\qed

The primitive elements (in the duplicial sense) are defined by the condition $\delta x=0$.
It is known that the binary operation
\begin{equation}\label{eq:mag}
\{x,y\} = x\prec y -x\succ y
\end{equation}
is magmatic, and that it preserves the primitive subalgebra 
of a duplicial algebra. As a consequence, the $Dup$-primitive
subalgebra of $\CQSym$ is generated by $x=\P^1$ for the operation \eqref{eq:mag}
This implies that $(As,Dup,Mag)$ is a good triple of operads in the
sense of Loday \cite{LodTrip}. That is, we have $As\circ Mag=Dup$.

The first primitive elements are
\begin{eqnarray}
\{x,x\} &=& \P^{11}-\P^{12},\\
\{\{x,x\},x\} &=& \P^{111}-\P^{122}-\P^{113}+\P^{123},\\
\{x,\{x,x\}\} &=&  \P^{111}-\P^{122}-\P^{112}+\P^{123}.
\end{eqnarray}

As an associative algebra for the product $\succ$, $\CQSym$ is free
over its primitive subspace (which has dimension $c_{n-1}$ in degree
$n$ as expected).

The duplicial operad is Koszul, and its dual $Dup^!$ is defined
by the same relations, together with two extra ones
\begin{equation}
(x\prec y)\succ z=0\quad\text{and}\quad 0=x\prec(y\succ z)
\end{equation}
The dimension of $Dup^!$  in degree $n$ is $n$, and a linear basis is formed by
the hook-shaped trees
\begin{equation}
x\succ x\succ \cdots \succ x \prec x\prec\cdots\prec x\prec x.
\end{equation}

There are known morphisms of operads  $Dup\rightarrow As^{(2)}$,
$Dup\rightarrow Dias$ and $2as\rightarrow Dup$ \cite{LodTrip}.
In particular, the morphism of operads $Dup\rightarrow As^{(2)}$ corresponds to
the Hopf algebra morphism $\phi:\ \CQSym\rightarrow \Sym$ 
(defined in \cite{NTPark})
\begin{equation}
\phi(\P^\pi)=S^{t(\pi)}
\end{equation}
where the composition $t(\pi)$ is the packed evaluation  $\pev(\pi)$, i.e.,
the composition obtained by removing the zeros from $\ev(\pi)$.

%%%%%%%%%%%%%%%%%%%%%%%%%%%%%%%%%%%%%%%%%%%%%%%%%%%%%%%%%%%%%%%%%%%%%%%%%%%%%%%
\subsection{Free quasi-symmetric functions and dendriform algebras}

A dendriform algebra~\cite{Lod0,Lod} is an associative algebra $A$ 
endowed with two bilinear operations $\gaudend$, $\droitdend$,
satisfying
\begin{eqnarray}
(x\gaudend y)\gaudend z &=& x\gaudend (y\cdot z)\,,\\
(x\droitdend y)\gaudend z &=& x\droitdend (y\gaudend z)\,,\\
(x\cdot y)\droitdend z &=& x\droitdend (y\droitdend z)\,,
\end{eqnarray}
such that the associative
multiplication $\cdot$ splits as
\begin{equation}
a\cdot b = a \gaudend b + a \droitdend b\,.
\end{equation}
If the algebra has a unit, one has to add the conditions
\begin{equation}
1\succ x=x,\quad 1\prec x = 0,\quad x\succ 1 = 0,\quad x\prec 1 = x\quad (x\not=1)
\end{equation}
and leave $1\prec 1$ and $1\succ 1$ undefined.

The algebra of free quasi-symmetric functions $\FQSym$~\cite{NCSF6} (or the
Malvenuto-Reutenauer Hopf algebra of permutations~\cite{MR}) is dendriform.

For a totally ordered alphabet $A$, $\FQSym(A)$ is the algebra spanned
by the noncommutative ``polynomials''\footnote{
If $A$ is infinite, these are elements of an inverse limit of noncommutative
polynomial algebras.
}
\begin{equation}
\G_\sigma(A)  := \sum_{\genfrac{}{}{0pt}{}{w\in A^n}{\std(w)=\sigma}} w
\end{equation}
where $\sigma$ is a permutation in the symmetric group $\SG_n$ and $\std(w)$
denotes the standardization of the word $w$.
The multiplication rule is 
\begin{equation}\label{eq:prodG}
\G_\alpha \G_\beta = \sum_{\gamma\in \alpha * \beta} \G_\gamma,
\end{equation}
where the convolution $\alpha*\beta$ 
of $\alpha\in\SG_k$ and $\beta\in\SG_l$
is the (multiplicity free) sum  in the group algebra of $\SG_{k+l}$~\cite{MR}
\begin{equation}
\alpha * \beta =
\sum_{\genfrac{}{}{0pt}{}{\gamma=u\cdot v}{\std(u)=\alpha ;\, \std(v)=\beta}}
\gamma\,,
\end{equation}
regarded as a set in \eqref{eq:prodG}.
The dendriform structure of $\FQSym$ is given by
\begin{equation}
\G_\alpha \G_\beta = \G_\alpha \gaudend \G_\beta + \G_\alpha \droitdend
\G_\beta\,,
\end{equation}
where
\begin{equation}
\G_\alpha \gaudend \G_\beta =
\sum_{\genfrac{}{}{0pt}{}{\gamma=u\cdot v \in \alpha*\beta}%
{|u|=|\alpha| ;\, \max(v)<\max(u)}}
\G_\gamma\,,
\end{equation}
\begin{equation}
\G_\alpha \droitdend \G_\beta =
\sum_{\genfrac{}{}{0pt}{}{\gamma=u.v\in \alpha*\beta}%
{|u|=|\alpha| ;\, \max(v)\geq\max(u)}} \G_\gamma\,.
\end{equation}
There is a scalar product on $\FQSym$ given by
\begin{equation}
\<\G_\sigma,\G_\tau\>=\delta_{\sigma,\tau^{-1}}
\end{equation}
so that
\begin{equation}
\F_\sigma := \G_{\sigma^{-1}}
\end{equation}
is the dual basis of $\G_\sigma$ (as a Hopf algebra, $\FQSym$
is self-dual).

Now, $x=\G_1$ generates a free dendriform algebra in $\FQSym$, 
$\PBT$, the Loday-Ronco algebra of planar binary trees~\cite{LR1}.
The natural basis $\P_T$ of this algebra can be interpreted as follows.

In $\FQSym$, one can build from the dendriform operations a bilinear map
\cite{HNT}
\begin{equation}\label{defB}
B(F,G) = F\succ \G_1 \prec G
\end{equation}
so that the terms $B_T(a)$ of the binary tree expansion of the unique solution
of the functional equation
\begin{equation}\label{FE}
X = a + B(X,X)
\end{equation}
are precisely the basis $\P_T$ of the free dendriform algebra on one generator
$x=\G_1$ for the choice $a=1$ ($B_T(a)$ is the result of evaluating the
expression encoded by the complete binary tree $T$ with $a$ in the leaves and
$B$ in the internal nodes).
As for the solution $X$, it is just the sum of all $\G_\sigma$, that is, the
sum of all words.

This interpretation leads to simple derivations of the $q$-hook length
formulas for binary trees \cite{HNT}, as well as to the combinatorial
interpretations of various special functions such as the tangent \cite{JVNT}
or possibly the Jacobi elliptic functions.

%%%%%%%%%%%%%%%%%%%%%%%%%%%%%%%%%%%%%%%%%%%%%%%%%%%%%%%%%%%%%%%%%%%%%%%%%%%%%%%
\subsection{Word quasi-symmetric functions and tridendrifom algebras}

A \emph{dendriform trialgebra}~\cite{LRtri}
(or tridendriform algebra)  is an associative algebra whose multiplication
$\odot$ splits into three pieces
\begin{equation}
x\odot y = x\gautrid y + x\miltrid y + x\droittrid y\,,
\end{equation}
where $\miltrid$ is associative, and
\begin{eqnarray}
(x\gautrid y)\gautrid z = x\gautrid (y\odot z)\,,\\
(x\droittrid y)\gautrid z = x\droittrid (y\gautrid z)\,,\\
(x\odot y)\droittrid z = x\droittrid (y\droittrid z)\,,\\
(x\droittrid y)\miltrid z = x\droittrid (y\miltrid z)\,,\\
(x\gautrid y)\miltrid z = x\miltrid (y\droittrid z)\,,\\
(x\miltrid y)\gautrid z = x\miltrid (y\gautrid z)\,.
\end{eqnarray}

The \emph{packed word} $u=\pack(w)$ associated with a word $w\in A^*$ is
obtained by the following process. If $b_1<b_2<\ldots <b_r$ are the letters
occuring in $w$, $u$ is the image of $w$ by the homomorphism
$b_i\mapsto a_i$.
A word $u$ is said to be \emph{packed} if $\pack(u)=u$. We denote by $\PW$ the
set of packed words.
With such a word, we associate the polynomial
\begin{equation}
\M_u :=\sum_{\pack(w)=u}w\,.
\end{equation}

These polynomials span a subalgebra of $\K\langle A\rangle$, called $\WQSym$
for Word Quasi-Symmetric functions see, {\it e .g.}, \cite{HNT}).

The product on $\WQSym$ is given by
\begin{equation} 
\label{prodG-wq}
\M_{u'} \M_{u''} = \sum_{u \in u'\convW u''} \M_u\,,
\end{equation}
where the \emph{convolution} $u'\convW u''$ of two packed words
is defined as
\begin{equation} 
u'\convW u'' = \sum_{v,w ;
u=v\cdot w\,\in\,\PW, \pack(v)=u', \pack(w)=u''} u\,.
\end{equation}
It is a dendriform trialgebra. 
The partial products are given by
\begin{equation}
\M_{w'} \gautrid \M_{w''} =
\sum_{w=u.v\in w'\convW w'', |u|=|w'| ; \max(v)<\max(u)}
\M_\park,
\end{equation}
\begin{equation}
\M_{w'} \miltrid \M_{w''} =
\sum_{w=u.v\in w'\convW w'', |u|=|w'| ; \max(v)=\max(u)}
\M_\park,
\end{equation}
\begin{equation}
\M_{w'} \droittrid \M_{w''} =
\sum_{w=u.v\in w'\convW w'', |u|=|w'| ; \max(v)>\max(u)}
\M_\park,
\end{equation}
It is known~\cite{LRtri} that the free dendriform trialgebra on one
generator, denoted here by
$\TD$, is a free associative algebra with Hilbert series
\begin{equation}
\label{sg-dendt}
\sum_{n\geq0} s_n t^n = \frac{1+t-\sqrt{1-6t+t^2}}{4t}
= 1 + t + 3t^2 + 11t^3 + 45t^4 + 197t^5 + \cdots
\end{equation}
that is, the generating function of the \emph{super-Catalan}, or
\emph{little Schr\"oder} numbers, counting \emph{reduced plane trees},
{\it i.e.,}, plane rooted trees in which internal nodes have at least two
descendants \cite[A001003]{Slo}.

There is a natural embedding $\FQSym\hookrightarrow\WQSym$
given by
\begin{equation}
\G_\sigma \mapsto \sum_{\std(u)=\sigma}\M_u\,.
\end{equation}
On the polynomial realizations, this is indeed an inclusion.

Both $\FQSym$ and $\WQSym$ can be interpreted as operads. The
space of $n$-ary operations of the Zinbiel operad  can be naturally
identified with $\FQSym_n$ \cite{CHNT}, and for $\WQSym$, the
relevant operad is described in \cite{MNT}.

%%%%%%%%%%%%%%%%%%%%%%%%%%%%%%%%%%%%%%%%%%%%%%%%%%%%%%%%%%%%%%%%%%%%%%%%%%%%%%%
\subsection{The triduplicial operad}

The free dendriform and free tridendriform algebras
on one generator arise naturally as subalgebras and quotients of
$\FQSym$ and  $\WQSym$, as well
as the free tricubical algebra. More precisely, the sylvester
quotient of $\WQSym$ is isomorphic to the free tridendriform algebra on one generator
(as a tridendriform algebra and as a Hopf algebra), and its hypoplactic quotient is
isomorphic to the free tricubical algebra on one generator \cite{NTWQS}.

 Applying the same strategy to
$\PQSym$, we obtain together with $\CQSym$ (Catalan numbers),
the little Schr\"oder numbers ($\SQSym$) and segmented compositions
(powers of 3) \cite{NTPark}. Having related $\CQSym$ to the duplicial operad,
we may suspect that $\SQSym$ should
have an operadic interpretation which is to $Dup$ what $TriDend$
is to $Dend$. This is indeed the case. The new operad will
be called \emph{triduplicial}. It is a quotient of the triplicial
operad defined in \cite{Ler} (one more relation).

\begin{definition}
\label{def-tridup}
A triduplicial algebra is a vector space $V$ endowed with three
\emph{associative} laws $\prec,\succ,\circ$ such that\\
(i) $(x\succ y)\prec z=x\succ(y\prec z)$, i.e., $(V,\prec,\succ)$ is duplicial,\\
%(i) $(V,\prec,\succ)$ is a duplicial algebra\\
(ii) $(x\circ y)\prec z=x\circ(y\prec z)$, i.e., $(V,\prec,\circ)$ is duplicial,\\
(iii) $(x\succ y)\circ z=x\succ (y\circ z)$, i.e., $(V,\circ,\succ)$ is duplicial,\\
(iv) $(x\circ y)\succ z=x\circ(y\succ z)$, i.e., $(V,\succ,\circ)$ is duplicial.
\end{definition}

Recall from \cite{NTPark} that a {\it parking quasi-ribbon} is 
a segmented nondecreasing parking function where the bars only occur at
positions $\cdots a\sep b\cdots$, with $a<b$. These objects encode
hypoplactic classes of parking functions. The first ones are
\begin{equation}
\label{ex12}
\{1\}, \qquad\qquad \{11,\, 12,\, 1\sep2\},
\end{equation} 
\begin{equation}
\label{ex3}
\{111,\  112,\  11\sep2, \  113,\  11\sep3,\  122,\  1\sep22,\ 
123, \ 1\sep23, \ 12\sep3, \ 1\sep2\sep3 \}.
\end{equation}
The number of parking quasi-ribbons of length $n$ is the little Schr\"oder
number $s_n$. 

With a parking quasi-ribbon $\qrpark$, we associate the elements
\begin{equation}
\PS_{\qrpark} := \sum_{{\sf P}(\park)=\qrpark} {\F_\park},
\end{equation}
where ${{\sf P}}(\a)$ denotes the hypoplactic class of $\a$.
For example,
\begin{equation}
\PS_{11|3} = \F_{131}+\F_{311}\,, \qquad
\PS_{113} = \F_{113}.
\end{equation}
The $\PS_\qrpark$ form a basis of a Hopf subalgebra of $\PQSym$,
denoted by $\SQSym$ \cite{NTPark}. As an associative algebra, it
is isomorphic to $\TD$, the free tridendriform algebra on one generator,
which is itself a Hopf subalgebra of $\WQSym$ constructed by a similar
method. However, this is not an isomorphism of Hopf algebras.
Indeed, $\TD$ is self-dual, but $\SQSym$ is not. This raised the question
of an operadic interpretation of $\SQSym$. The following result
provides an answer.
 
\begin{theorem}
The free triduplicial algebra on one generator $\T$  has the little Schr\"oder
numbers as graded dimensions. Its natural basis can be realised by
parking quasi-ribbons,
if we define $\prec$ and $\succ$ by shifted concatenation as above, and
$\circ$ as ordinary shifted concatenation with insertion of a bar:
\begin{equation}
\P_{\q'}\circ \P_{\q''}= \P_{\q'|\q''[|\q'|]}\,.
\end{equation}
\end{theorem}

\Proof
It is immediate to check that this defines a triduplicial structure, and that
the subspace of $\SQSym$ generated by $\P_1$ for these three
operations contains all the quasi-ribbons $\P_\q$.
So the free triduplicial algebra has at least the sequence
$s_n$ as graded dimensions.

Now, the triduplicial relations can be presented as rewriting rules for evaluation
trees. We then have the seven relations, the first three consisting in the
associativity of the three operations, the next four ones being the
 four duplicial relations presented in
Definition~\ref{def-tridup}.

\newdimen\vcadre\vcadre=0.01cm % marges verticales de la boite
\newdimen\hcadre\hcadre=0.01cm % marges horizontales de la boite
\begin{equation}
%\small
\arbg{<}{<} \longleftrightarrow \arbd{<}{<}
\qquad\quad
\arbg{\circ}{\circ} \longleftrightarrow \arbd{\circ}{\circ}
\end{equation}

\begin{equation}
\arbg{>}{>} \longleftrightarrow \arbd{>}{>}
\end{equation}

\begin{equation}
\arbg{<}{\circ} \longleftrightarrow \arbd{\circ}{<}
\qquad\quad
\arbg{<}{>} \longleftrightarrow \arbd{>}{<}
\end{equation}

\begin{equation}
\arbg{\circ}{>} \longleftrightarrow \arbd{>}{\circ}
\qquad\quad
\arbg{>}{\circ} \longleftrightarrow \arbd{\circ}{>}
\end{equation}

Now consider all these relations as oriented rewriting between trees, each left
tree being replaced by its right counterpart. Then the free triduplicial algebra
on one generator is spanned (not necessarily freely)
by the trees that cannot be rewritten.
Analysing the seven relations above, one sees that the trees that cannot be
rewritten are of the following form:
either the tree consisting in a single node, or a tree with  any operation at
its root whose left subtree is a leaf, or one of the following two trees:
\begin{equation}
\arbgnew{\circ}{<}
\qquad\qquad
\arbgnew{>}{<},
\end{equation}
where $A$ and $B$ are trees that cannot be rewritten.
The generating series $S(x)$ of these trees satisfies therefore
the functional equation
\begin{equation}
S = 1+ 3xS + 2x^2S^2,
\end{equation}
whose solution is  the generating series of the $s_n$, so that that
the free triduplicial algebra on one generator has at most the $s_n$ as graded
dimensions.  Therefore, the parking quasi-ribbons provide a faithful realization of
it.
\qed

Thus, $\T$ can be identified with $\SQSym$. Another description of $\SQSym$
can be found in \cite{MNTconj}, in terms of {\it Schr\"oder pseudocompositions}
(the Polish codes of reduced plane trees). In this reference, Poincar\'e's functional
equation is cast into the form of a $q$-deformation of \eqref{FEG} below, allowing
to state that $\SQSym$ is to Poincar\'e's equation what $\CQSym$ is to Lagrange
inversion.

For the sake of completeness, let us mention that as in the case of $Dup$,
$TriDup$ is Koszul, and the
dual operad $TriDup^!$ is defined by the same relations, together with four
extra ones\footnote{
These facts have been proved in 2012 by A. Mansuy (unpublished).
}
\begin{eqnarray}
(x\prec y)\succ z=0\quad&\text{and}&\quad 0=x\prec(y\succ z),\\
(x\prec y)\circ z=0\quad&\text{and}&\quad 0=x\prec(y\circ z).
\end{eqnarray}

The sequence of dimensions of $TriDup^!$ is therefore $2^n-1$, as for the
dual of {TriDend} \cite{Zin}.

%%%%%%%%%%%%%%%%%%%%%%%%%%%%%%%%%%%%%%%%%%%%%%%%%%%%%%%%%%%%%%%%%%%%%%%%%%%%%%
\subsection{Duplicial operations on parking functions}

The free dendriform and tridendriform algebras on one
generator arise naturally as subalgebras of $\FQSym$
and $\WQSym$, which are themselves  dendriform and tridendriform
in a natural way.
Similarly, the duplicial structure of the Catalan algebra
$\CQSym$ is actually inherited from
a duplicial structure on $\PQSym$. The $\succ$ operation
is the usual product (given by the ordinary shifted shuffle).
The $\prec$ operation is also a kind of shifted shuffle,
with a normalization factor. 
\begin{proposition}
If $\max(\a)=m$, let
\begin{equation}
\F_\a\prec\F_\b
= \frac{|\a|_m!|\b|_1!}{(|\a|_m+|\b|_1)!}\F_{\a\shuffle \b[m-1]}\,.
\end{equation}
Then, $\PQSym$ is duplicial for $\prec$ and $\succ$.
\end{proposition}

\Proof Direct verification. \qed

%%%%%%%%%%%%%%%%%%%%%%%%%%%%%%%%%%%%%%%%%%%%%%%%%%%%%%%%%%%%%%%%%%%%%%%%%%%%%%%
%%%%%%%%%%%%%%%%%%%%%%%%%%%%%%%%%%%%%%%%%%%%%%%%%%%%%%%%%%%%%%%%%%%%%%%%%%%%%%%
%%%%%%%%%%%%%%%%%%%%%%%%%%%%%%%%%%%%%%%%%%%%%%%%%%%%%%%%%%%%%%%%%%%%%%%%%%%%%%%
\section{Lagrange inversion}

%%%%%%%%%%%%%%%%%%%%%%%%%%%%%%%%%%%%%%%%%%%%%%%%%%%%%%%%%%%%%%%%%%%%%%%%%%%%%%%
\subsection{A bilinear duplicial equation} \label{newbij}
Let $G\in\PQSym$ be the formal sum of all parking functions 
\begin{equation}
G = \sum_{\a\in\PF}\F_\a\,.
\end{equation}
Actually, $G$ belongs to $\CQSym$, and
\begin{equation}
G = \sum_{\pi\in\NDPF}\P^\pi\,.
\end{equation}

\begin{proposition}
Define a bilinear map $B$ on $\CQSym$ by formula (\ref{defB}),  interpreting
now $\prec$ and $\succ$ as the duplicial operations. Then, $G$ satisfies
the functional equation
\begin{equation}\label{FEG}
G = 1 + B(G,G)
\end{equation}
and each term $B_T(1)$ of the tree expansion of the solution is
a single $\P^\pi$, thus forcing a bijection beween binary trees
and nondecreasing parking functions.
\end{proposition}

\Proof
Proposition~\ref{prop-cqs} shows that $\CQSym$ is the free duplicial algebra
on one generator. In particular,  each nondecreasing parking
function $\pi$  has a unique expression of the form $\P^\alpha\succ\P^1\prec
\P^\beta$, where $\alpha$ and $\beta$ are (possibly empty) nondecreasing
parking functions.

The same is true of binary trees, if one interprets $\prec$ and $\succ$ as the
over-under operations~\cite{LR1,LodTrip}, whence the correspondence.
\qed

The bijection between binary trees and nondecreasing parking 
functions can be described as follows. Starting with a tree $T$, its vertices
are recursively labeled by integers, and the tree is flattened so as to read
a word. The label of the root is the number $m$ of vertices of its left
subtree, plus one. The labels of the right subtree are its original ones
shifted by $m-1$.

For example,
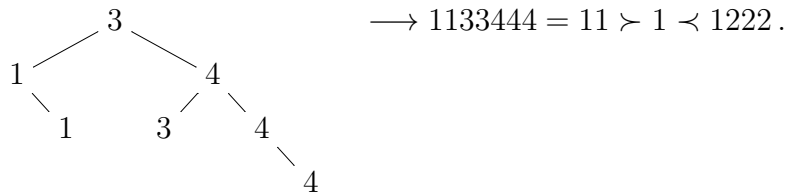
\begin{equation*}
{\xymatrix@C=2mm@R=2mm{
        *{} & *{} & {3}\ar@{-}[drr]\ar@{-}[dll] \\
        {1}\ar@{-}[dr]  & *{} & *{} & *{} & {4}\ar@{-}[dr]\ar@{-}[dl] \\
        *{}  & {1} & *{} & {3} & *{} & 4\ar@{-}[dr] \\
        *{}  & *{} & *{} & *{} & *{} & *{} & 4
      }}
\quad \longrightarrow 1133444 = 11\succ 1\prec 1222\,.
\end{equation*}

%%%%%%%%%%%%%%%%%%%%%%%%%%%%%%%%%%%%%%%%%%%%%%%%%%%%%%%%%%%%%%%%%%%%%%%%%%%%%%%
\subsection{Some Tamari intervals}

This correspondence has an interesting compatibility with the Tamari
order: 

\begin{proposition}
Nondecreasing parking functions with the same packed evaluation $I$
form an interval, whose cardinality is the coefficient of $S^I$ in $g$.
\end{proposition}

\Proof It is enough to show that the nondecreasing parking functions
with packed evaluation $J\le I$ for the reverse refinement order form
an interval. Now, $1^{i_1}\succ\ldots\succ 1^{i_r}$ is maximal among those
for the Tamari order (this is the maximal element of the product in $\PBT$
of the trees encoded by the $1^{i_k}$), and
$1^n=1^{i_1}\prec\ldots\prec 1^{i_r}$ is its minimal element. This corresponds
to $S^I$ through the embedding of $\Sym$ given by $S_n\mapsto 1^n$. Thus
this is an interval, and the intervals composed of nondecreasing parking
functions with the same packed evaluation $I$ correspond to the expansion of
the ribbons $R_{\bar I^\sim}$ (conjugate mirror of $I$) on the basis of trees.
\qed

In other words, our bijection between binary trees and nondecreasing parking
functions has the property that the trees having the same canopy~\cite{Vie,LR1} 
correspond 
to  nondecreasing parking functions with the same packed evaluation.

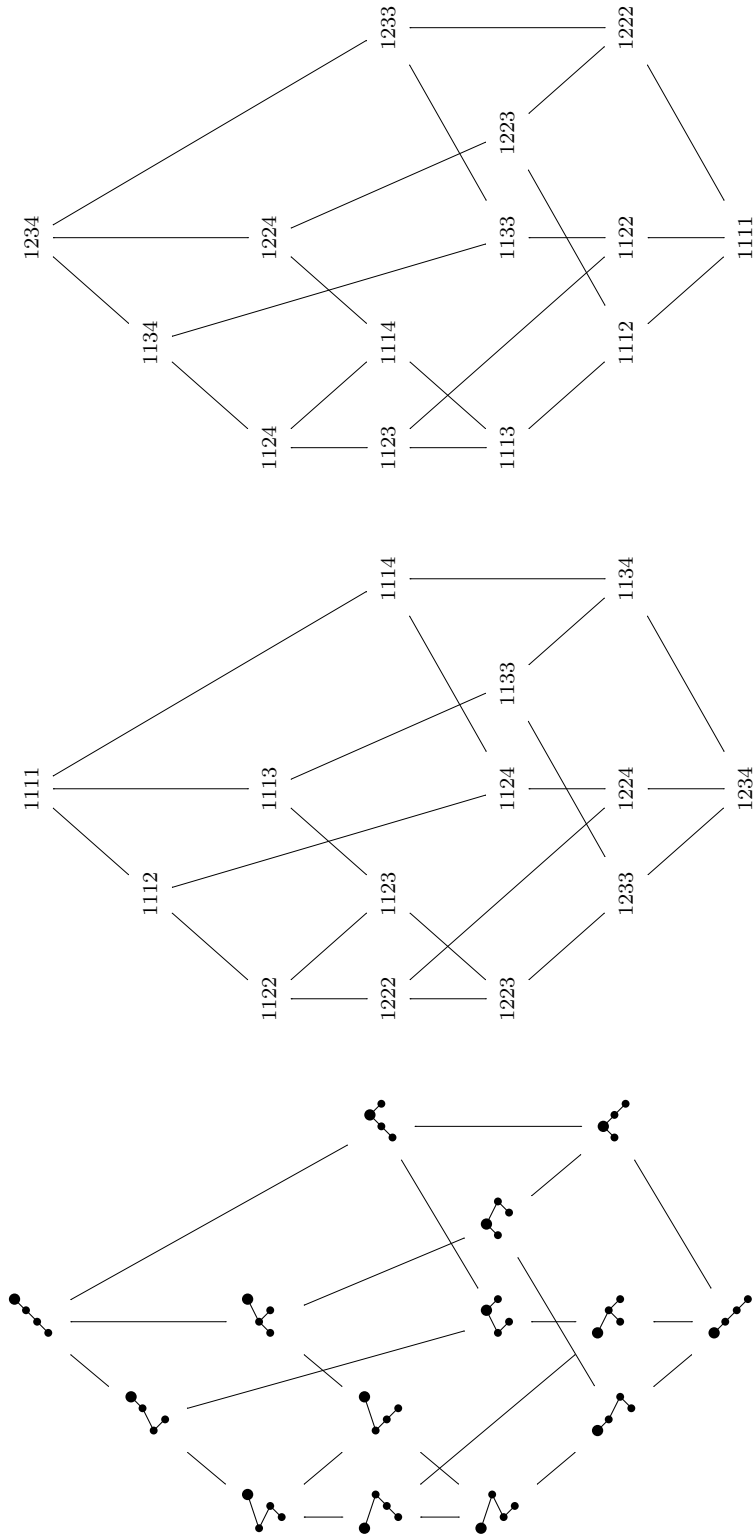
\begin{figure}[ht]
{\centerline{
\tiny
\rotateleft{
\newdimen\vcadre\vcadre=0.15cm % marges verticales de la boite
\newdimen\hcadre\hcadre=0.15cm % marges horizontales de la boite
\setlength\unitlength{1.5mm}
$\xymatrix@R=0.5cm@C=4mm{
 &  & *{\GrTeXBox{\begin{picture}(4,5)\put(1,1){\circle*{0.7}}\put(2,2){\circle*{0.7}}\put(2,2){\Line(-1,-1)}\put(3,3){\circle*{0.7}}\put(3,3){\Line(-1,-1)}\put(4,4){\circle*{0.7}}\put(4,4){\Line(-1,-1)}\put(4,4){\circle*{1}}\end{picture}}}\arx1[ld]\arx1[dd]\arx1[rrddd]& \\
 & *{\GrTeXBox{\begin{picture}(4,5)\put(1,2){\circle*{0.7}}\put(2,1){\circle*{0.7}}\put(1,2){\Line(1,-1)}\put(3,3){\circle*{0.7}}\put(3,3){\Line(-2,-1)}\put(4,4){\circle*{0.7}}\put(4,4){\Line(-1,-1)}\put(4,4){\circle*{1}}\end{picture}}}\arx1[ld]\arx1[rddd]& \\
*{\GrTeXBox{\begin{picture}(4,5)\put(1,3){\circle*{0.7}}\put(2,1){\circle*{0.7}}\put(3,2){\circle*{0.7}}\put(3,2){\Line(-1,-1)}\put(1,3){\Line(2,-1)}\put(4,4){\circle*{0.7}}\put(4,4){\Line(-3,-1)}\put(4,4){\circle*{1}}\end{picture}}}\arx1[d]\arx1[rd]&  & *{\GrTeXBox{\begin{picture}(4,5)\put(1,2){\circle*{0.7}}\put(2,3){\circle*{0.7}}\put(3,2){\circle*{0.7}}\put(2,3){\Line(-1,-1)}\put(2,3){\Line(1,-1)}\put(4,4){\circle*{0.7}}\put(4,4){\Line(-2,-1)}\put(4,4){\circle*{1}}\end{picture}}}\arx1[ld]\arx1[rdd]& \\
*{\GrTeXBox{\begin{picture}(4,5)\put(1,4){\circle*{0.7}}\put(2,1){\circle*{0.7}}\put(3,2){\circle*{0.7}}\put(3,2){\Line(-1,-1)}\put(4,3){\circle*{0.7}}\put(4,3){\Line(-1,-1)}\put(1,4){\Line(3,-1)}\put(1,4){\circle*{1}}\end{picture}}}\arx1[d]\arx1[rrdd]& *{\GrTeXBox{\begin{picture}(4,5)\put(1,3){\circle*{0.7}}\put(2,2){\circle*{0.7}}\put(3,1){\circle*{0.7}}\put(2,2){\Line(1,-1)}\put(1,3){\Line(1,-1)}\put(4,4){\circle*{0.7}}\put(4,4){\Line(-3,-1)}\put(4,4){\circle*{1}}\end{picture}}}\arx1[ld]&  &  & *{\GrTeXBox{\begin{picture}(4,4)\put(1,1){\circle*{0.7}}\put(2,2){\circle*{0.7}}\put(2,2){\Line(-1,-1)}\put(3,3){\circle*{0.7}}\put(4,2){\circle*{0.7}}\put(3,3){\Line(-1,-1)}\put(3,3){\Line(1,-1)}\put(3,3){\circle*{1}}\end{picture}}}\arx1[lld]\arx1[dd]& \\
*{\GrTeXBox{\begin{picture}(4,5)\put(1,4){\circle*{0.7}}\put(2,2){\circle*{0.7}}\put(3,1){\circle*{0.7}}\put(2,2){\Line(1,-1)}\put(4,3){\circle*{0.7}}\put(4,3){\Line(-2,-1)}\put(1,4){\Line(3,-1)}\put(1,4){\circle*{1}}\end{picture}}}\arx1[rd]&  & *{\GrTeXBox{\begin{picture}(4,4)\put(1,2){\circle*{0.7}}\put(2,1){\circle*{0.7}}\put(1,2){\Line(1,-1)}\put(3,3){\circle*{0.7}}\put(4,2){\circle*{0.7}}\put(3,3){\Line(-2,-1)}\put(3,3){\Line(1,-1)}\put(3,3){\circle*{1}}\end{picture}}}\arx1[d]& *{\GrTeXBox{\begin{picture}(4,4)\put(1,2){\circle*{0.7}}\put(2,3){\circle*{0.7}}\put(3,1){\circle*{0.7}}\put(4,2){\circle*{0.7}}\put(4,2){\Line(-1,-1)}\put(2,3){\Line(-1,-1)}\put(2,3){\Line(2,-1)}\put(2,3){\circle*{1}}\end{picture}}}\arx1[lld]\arx1[rd]& \\
 & *{\GrTeXBox{\begin{picture}(4,5)\put(1,4){\circle*{0.7}}\put(2,3){\circle*{0.7}}\put(3,1){\circle*{0.7}}\put(4,2){\circle*{0.7}}\put(4,2){\Line(-1,-1)}\put(2,3){\Line(2,-1)}\put(1,4){\Line(1,-1)}\put(1,4){\circle*{1}}\end{picture}}}\arx1[rd]& *{\GrTeXBox{\begin{picture}(4,5)\put(1,4){\circle*{0.7}}\put(2,2){\circle*{0.7}}\put(3,3){\circle*{0.7}}\put(4,2){\circle*{0.7}}\put(3,3){\Line(-1,-1)}\put(3,3){\Line(1,-1)}\put(1,4){\Line(2,-1)}\put(1,4){\circle*{1}}\end{picture}}}\arx1[d]&  & *{\GrTeXBox{\begin{picture}(4,4)\put(1,2){\circle*{0.7}}\put(2,3){\circle*{0.7}}\put(3,2){\circle*{0.7}}\put(4,1){\circle*{0.7}}\put(3,2){\Line(1,-1)}\put(2,3){\Line(-1,-1)}\put(2,3){\Line(1,-1)}\put(2,3){\circle*{1}}\end{picture}}}\arx1[lld]& \\
 &  & *{\GrTeXBox{\begin{picture}(4,5)\put(1,4){\circle*{0.7}}\put(2,3){\circle*{0.7}}\put(3,2){\circle*{0.7}}\put(4,1){\circle*{0.7}}\put(3,2){\Line(1,-1)}\put(2,3){\Line(1,-1)}\put(1,4){\Line(1,-1)}\put(1,4){\circle*{1}}\end{picture}}}& \\
}$
\hfill
\newdimen\vcadre\vcadre=0.2cm % marges verticales de la boite
\newdimen\hcadre\hcadre=0.2cm % marges horizontales de la boite
$\xymatrix@R=1.0cm@C=4mm{
 &  & *{\GrTeXBox{1111}}\arx1[ld]\arx1[dd]\arx1[rrddd]& \\
 & *{\GrTeXBox{1112}}\arx1[ld]\arx1[rddd]& \\
*{\GrTeXBox{1122}}\arx1[d]\arx1[rd]&  & *{\GrTeXBox{1113}}\arx1[ld]\arx1[rdd]& \\
*{\GrTeXBox{1222}}\arx1[d]\arx1[rrdd]& *{\GrTeXBox{1123}}\arx1[ld]
  &  &  & *{\GrTeXBox{1114}}\arx1[lld]\arx1[dd]& \\
*{\GrTeXBox{1223}}\arx1[rd]&  & *{\GrTeXBox{1124}}\arx1[d]& *{\GrTeXBox{1133}}\arx1[lld]\arx1[rd]& \\
 & *{\GrTeXBox{1233}}\arx1[rd]& *{\GrTeXBox{1224}}\arx1[d]&  & *{\GrTeXBox{1134}}\arx1[lld]& \\
 &  & *{\GrTeXBox{1234}}& \\
}$
\hfill
\newdimen\vcadre\vcadre=0.2cm % marges verticales de la boite
\newdimen\hcadre\hcadre=0.2cm % marges horizontales de la boite
$\xymatrix@R=1.0cm@C=4mm{
 &  & *{\GrTeXBox{1234}}\arx1[ld]\arx1[dd]\arx1[rrddd]& \\
 & *{\GrTeXBox{1134}}\arx1[ld]\arx1[rddd]& \\
*{\GrTeXBox{1124}}\arx1[d]\arx1[rd]&  & *{\GrTeXBox{1224}}\arx1[ld]\arx1[rdd]& \\
*{\GrTeXBox{1123}}\arx1[d]\arx1[rrdd]& *{\GrTeXBox{1114}}\arx1[ld]
&  &  & *{\GrTeXBox{1233}}\arx1[lld]\arx1[dd]& \\
*{\GrTeXBox{1113}}\arx1[rd]&  & *{\GrTeXBox{1133}}\arx1[d]& *{\GrTeXBox{1223}}\arx1[lld]\arx1[rd]& \\
 & *{\GrTeXBox{1112}}\arx1[rd]& *{\GrTeXBox{1122}}\arx1[d]&  & *{\GrTeXBox{1222}}\arx1[lld]& \\
 &  & *{\GrTeXBox{1111}}& \\
}$
}
}}
\caption{\label{tam4}The Tamari order on trees,  on
nondecreasing parking functions regarded as noncrossing partitions,  and via the new
bijection.}
\end{figure}

For example, the coefficients of
\begin{equation}
g_4 = S^4+3S^{31}+2S^{22}+S^{13}+3S^{211}+2S^{121}+S^{112}+S^{1111}
\end{equation}
can be read on Figure~\ref{tam4}, where one can also easily check the interval
property.

\FloatBarrier
Each interval consists in the nonequivalent bracketings of a duplicial product
involving the operations $\prec$ and $\succ$ in the same order. For example,
the interval $1123-1124-1134$ consists in the bracketings of the word
$1 \prec 1 \succ 1 \succ 1$, which are
$1 \prec ((1 \succ 1) \succ 1)=1123$,
$(1 \prec (1 \succ 1)) \succ 1=1124$, and
$((1 \prec 1) \succ 1) \succ 1=1124$.
The cover relation consists in appliying one associativity relation. The other
edges of the Hasse diagram are obtained by changing one $\succ$ into a
$\prec$. Under the bijection with trees, both operations correspond to a
rotation:

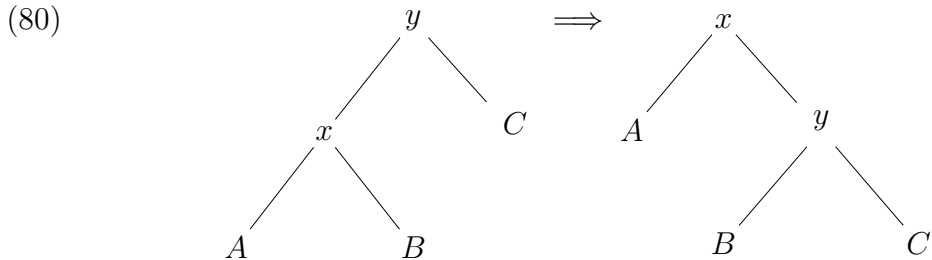
\begin{equation}
\newdimen\vcadre\vcadre=0.1cm % marges verticales de la boite
\newdimen\hcadre\hcadre=0.1cm % marges horizontales de la boite
\xymatrix@R=0.9cm@C=7mm{
 && *{\GrTeXBox{y}}\arx1[ld]\arx1[rd]\\
& *{\GrTeXBox{x}}\arx1[dl]\arx1[dr] && {\GrTeXBox{C}}\\
 *{\GrTeXBox{A}} &&  *{\GrTeXBox{B}}\\
}
\Longrightarrow
\xymatrix@R=0.9cm@C=7mm{
 & *{\GrTeXBox{x}}\arx1[ld]\arx1[rd]\\
*{\GrTeXBox{A}} & & {\GrTeXBox{y}}\arx1[dl]\arx1[dr]\\
& *{\GrTeXBox{B}}&  & *{\GrTeXBox{C}}\\
}
\end{equation}
where $A$, $B$, and $C$ are subtrees, 
the difference between both cases being whether the subtree $B$
is or not empty.
For example, $1113=1\prec((1\prec 1)\succ 1)$
and $1123=1\prec((1\succ 1)\succ 1)$.

%%%%%%%%%%%%%%%

\subsection{An application: counting sylvester classes}
Recall that the sylvester congruence  $\equiv$ is defined by
\begin{equation}
zxuy  \equiv xzuy\,,\ x\le y< z\in A\,,\ u\in A^*\,.
\end{equation}
The sylvester class of a word $w$ is a binary search tree $\PS(w)$,
which can be identified with a naked binary tree $T$ when $w$ is
a permutation. 

The natural basis $\P_T$ of $\PBT$, regarded as a subalgebra of $\FQSym$, can be described as \cite{HNT2}
\begin{equation}
\P_T= \sum_{\PS(\sigma)=T}\F_\sigma
\end{equation}
Let ${\rm can}(T)$ be the canonical permutation of the sylvester class indexed by $T$,
that is, the smallest permutation in the lexicographic order such that $\PS(\sigma)=T$.
The dual basis $\Q_T=\P_T^*$ can be realized as
\begin{equation}
\Q_T = \overline{\G_\sigma} =\G_\sigma(\bar A)
=\sum_{\std(u)=\sigma}\M_u(\bar A),
\end{equation}
where $\sigma={\rm can}(T)$, and
$\bar A$ is the image of our underlying alphabet by the canonical projection
$A^*\rightarrow A^*/\equiv$, 
so that its projection to $QSym$ (by the dual of the natural inclusion $\Sym\hookrightarrow\PBT$)
is
\begin{equation}
Q_T(X) = 
\sum_{I\in E(T)} M_I(X)
\end{equation}
where $E(T)$ is the set of evaluations of packed words $u$ such that $\std(u)={\rm can}(T)$.
Thus,
\begin{equation}
\sum_I\nu_IM_I(X)=\sum_T Q_T(X)
\end{equation}
and denoting by $\langle\cdot,\cdot\rangle$ the duality bracket between $QSym$ and $\Sym$,
 the number of sylvester classes of words of evaluation $I$ is
\begin{align}
\nu_I &= \sum_T \<Q_T,S^I\> \ \text{(sum over all binary trees)}\\
&= \sum_T\sum_{J\le I}\<Q_T,R_J\>\\
&= \sum_{J\le I}\sharp\NDPF(\bar J^\sim) \quad\text{thanks to the bijection of Section \ref{newbij}}\\
&=   \sum_{J\le I}\<M_{\bar J^\sim},g\>\quad \text{(by \eqref{gensol})}.
\end{align}
Since the coefficient of $S^{\bar J}$ in $g$ is equal to that of $S^{\bar J^\sim}$ \cite{NTLag}, see also
Section \ref{sec:invol},
we have finally:

\begin{proposition}
The number of sylvester classes of evaluation $I$ is
\begin{equation}
\nu_I =  \sum_{J\le I}\<M_{\bar J^\sim},g\> =  \sum_{J\le I}\<M_{\bar J},g\>.
\end{equation}
\end{proposition}

Defining a basis $E_I$ of $QSym$ by
\begin{equation}
E_I=\sum_{J \le I}M_J
\end{equation}
and denoting by $L_I$ its dual basis in $\Sym$, we can rewrite this number as
\begin{equation}\label{eq:nbslv}
\<E_{\overline{I}},g\>\ =\ \text{coefficient of $L_{\overline{I}}$ in $g$}.
\end{equation}

For example, the number of sylvester classes of permutations is
\begin{equation}
\nu_{1^n} = \<E_{1^n},g\> = \sum_{J\vDash n}\<M_I,g\>
\end{equation}
which is the sum of the coefficients of $g_n$ in the basis
$S^I$, hence the $n$th Catalan number. 

To find the number of sylvester classes of packed words,
we have to compute
\begin{equation}
\left\< \sum_{I\vDash n}E_{\overline{I}},g\right\>=
\left\< \sum_{\lambda\vdash n}2^{n-\ell(\lambda)}m_\lambda,g\right\>
\end{equation}
where the $m_\lambda$ are the monomial symmetric functions.
This can be evaluated by noting that
\begin{equation}
\begin{split}
\sum_{\lambda\vdash n}2^{n-\ell(\lambda)}m_\lambda
&= \frac1{(1-t)^n}
   \sum_{\lambda\vdash n}(1-t)^{\ell(\lambda)}m_\lambda|_{t=1/2} \\
&=\frac1{(1-t)^n} h_n((1-t)X)|_{t=1/2}
\end{split}
\end{equation}
Thus, the number of sylvester classes of packed words
is obtained by putting $x=2$ in the polynomial
\begin{equation}
\begin{split}
N_n(x)
&= \left. \frac1{(1-t)^n}\frac{h_n((1-t)(n+1))}{n+1}\right|_{t=1-1/x} \\
&= \frac1{n+1}
   \sum_{k=0}^n \binom{n+1}{k}\binom{2n-k}{n-k} (1-x)^kx^{n-k}
\end{split}
\end{equation}
a (reciprocal) Narayana polynomial, which gives back the little Schr\"oder
numbers, as expected.

%%%%%%%%%%%%%%%%%%%%%%%%%%%%%%%%%%%%%%%%%%%%%%%%%%%%%%%%%%%%%%%%%%%%%%%%%%%%%%%
\subsection{An involution on $\CQSym$}\label{sec:invol}

The natural involution (mirror symmetry) on binary trees
induces an involution $\iota$ on nondecreasing parking functions,
which is similar to (but different from) the involution $\nu$
of \cite{NTLag}.

Defining a basis $\Q^\pi$ by
\begin{equation}
\Q^\pi =\iota(\P^\pi)
\end{equation}
we have
\begin{equation}
\Q^{\pi'}\Q^{\pi''}=\Q^{\pi'\prec\pi''}\,.
\end{equation}

Now, there is a Hopf algebra morphism $\phi:\ \CQSym\rightarrow\Sym$
mapping $\P^\pi$ to $S^{t(\pi)}$, where the composition $t(\pi)$
is the packed evaluation of $\pi$.
This maps sends $G$ to the noncommutative Lagrange series $g$ \cite{NTLag},

\begin{equation}
\phi(G) = g\,,\ \text{the unique solution of}\ g=\sum_{n\ge 0}S_n g^n\,. 
\end{equation}

The involution $\iota$ is mapped to $S^I\mapsto S^{I^\sim}$. This provides
a simple proof of the symmetry of $g$ observed in \cite[Sec. 7]{NTLag}.
Here are the first values of $\iota$:

%\FloatBarrier
%\begin{figure}[ht]
\begin{equation}
\begin{tabular}{|c|c|}
\hline
$\pi$ & $\iota(\pi)$\\
\hline
1 & 1\\
\hline
12 & 11\\
\hline
123 & 111\\
113 & 112\\
122 & 122\\
\hline
\end{tabular}
\qquad
\begin{tabular}{|c|c|}
\hline
$\pi$ & $\iota(\pi)$\\
\hline
1234 & 1111\\
1134 & 1112\\
1224 & 1122\\
1124 & 1113\\
1114 & 1123\\
1233 & 1222\\
1133 & 1223\\
\hline
\end{tabular}
\end{equation}
%\caption{The first values of $\iota$.}
%\end{figure}
%%%%%%%%%%%%%%%%%%%%%%%%%%%%%%%%%%%%%%%%%%%%%%%%%%%%%%%%%%%%%%%%%%%%%%%%%%%%%%%
%%%%%%%%%%%%%%%%%%%%%%%%%%%%%%%%%%%%%%%%%%%%%%%%%%%%%%%%%%%%%%%%%%%%%%%%%%%%%%%
%%%%%%%%%%%%%%%%%%%%%%%%%%%%%%%%%%%%%%%%%%%%%%%%%%%%%%%%%%%%%%%%%%%%%%%%%%%%%%%
\section{Characters}

The identification the the Lagrange series as a homomorphic image of the sum of
all parking functions implies in particular that any combinatorial formula derived
by means of Lagrange inversion must come from a character of the algebra of 
parking functions. In this section, we illustrate this idea on a couple of examples.
%%%%%%%%%%%%%%%%%%%%%%%%%%%%%%%%%%%%%%%%%%%%%%%%%%%%%%%%%%%%%%%%%%%%%%%%%%%%%%%
\subsection{The formula of Lassalle for the Narayana polynomials}

It has been observed by M. Lassalle \cite{Lass} 
that the Narayana polynomials $c_n(q)$ could be expressed
in the form
\begin{equation}
qc_n(q)=\frac1{n+1}h_n((n+1)q)
\end{equation}
in $\lambda$-ring notation, with the assumption that $x=1-q$
is of rank 1.

Otherwise said,
\begin{equation}
(1-x)c_n(1-x)=g_n(1-x),
\end{equation}
the image of the symmetric function $g_n(X)$ by the character of $Sym$
sending the power sum $p_n$ to $1-x^n$. This character is also well-defined
on $\Sym$ \cite{NCSF2}, and we shall see that it can be lifted to $\PQSym$,
as well as the two-parameter extension
\begin{equation}
p_n\mapsto p_n\left(\frac{1-x}{1-t}\right)=\frac{1-x^n}{1-t^n}\,.
\end{equation}

This result can then be given a combinatorial
interpretation by following a chain of morphisms whose
composition builds up a character of $\PQSym$.

%%%%%%%%%%%%%%%%%%%%%%%%%%%%%%%%%%%%%%%%%%%%%%%%%%%%%%%%%%%%%%%%%%%%%%%%%%%%%%%
\subsection{From parking functions to permutations}

The polynomial realization of the dual $\PQSym^*$ given in  \cite{NT1}
implies the existence of an embedding of Hopf algebras
\begin{equation}
i:\ \G_\sigma\longmapsto \sum_{\std(\a)=\sigma}\G_\a
\end{equation}
which is actually an inclusion of the polynomial realizations:
\begin{equation}
\G_\sigma(A)= \sum_{\std(\a)=\sigma}\G_\a(A)\,.
\end{equation}
Thus, the dual map
\begin{equation}
i^*:\ \F_\a \longmapsto  \F_{\std(\a)}
\end{equation}
is a surjective morphism of Hopf algebras.

%%%%%%%%%%%%%%%%%%%%%%%%%%%%%%%%%%%%%%%%%%%%%%%%%%%%%%%%%%%%%%%%%%%%%%%%%%%%%%%
\subsection{From permutations to compositions}

The restriction of $i^*$ to $\SQSym$ takes its values in $\Sym$,
and precisely,
\begin{equation}
i^*(\P_\q) = R_I
\end{equation}
where $I$ is the shape of the quasi-ribbon $\q$.

Further restricting to $\CQSym\subset\SQSym$ yields
\begin{equation}
i^*(\P^\pi) = S^I
\end{equation}
where $I=t(\pi)$ is the packed evaluation of $\pi$.

%%%%%%%%%%%%%%%%%%%%%%%%%%%%%%%%%%%%%%%%%%%%%%%%%%%%%%%%%%%%%%%%%%%%%%%%%%%%%%%
\subsection{From compositions to scalars or polynomials}

At the level of  $\Sym$, we have many characters at our disposal.
Since characters take their values in a commutative algebra
(our ground field $\KK$), they come actually from characters
of $Sym$. An important and classical example is the evaluation
on the virtual alphabet
\begin{equation} 
\A = \frac{1-x}{1-q}\,\ \text{given by}\ p_n(\A)= \frac{1-x^n}{1-q^n}\,.
\end{equation}

This character can be lifted in several ways to $\FQSym$. One of
them, denoted by
\begin{equation} 
\A = \frac{\mid 1-x}{1-q\mid}
\end{equation}
is described in \cite{NTsuper} (see \eqref{qtFsigma} below).

%%%%%%%%%%%%%%%%%%%%%%%%%%%%%%%%%%%%%%%%%%%%%%%%%%%%%%%%%%%%%%%%%%%%%%%%%%%%%%%
\subsection{Super-Narayana polynomials}

Recall from \cite{NTcolored} that a signed parking function is a pair
$(\a,\varepsilon)$ where $\a$ is an ordinary parking function
and $\varepsilon$ a word of the same length over the alphabet $\{\pm 1\}$.
There are therefore $2^n(n+1)^{n-1}$ signed parking functions of
length $n$ (cf. \cite[A097629]{Slo}). There is a Hopf algebra,
denoted by $\PQSym^{(2)}$ based on signed parking functions. The product
in its $\F_{(\a,\varepsilon)}$ basis is given by the shifted shuffle
of signed words.

\begin{definition}
A signed inversion of a signed parking function $(\a,\varepsilon)$
is a pair $(i,j)$ with $i<j$ such that either
$\varepsilon_ia_i>\varepsilon_ja_j$ (an ordinary inversion) or
$\varepsilon_ia_i=\varepsilon_ja_j$ and the common sign $\varepsilon_i$ is
$-1$.

Similarly, a signed descent is an index $i$ such that either
$\varepsilon_ia_i>\varepsilon_{i+1}a_{i+1}$
or  $\varepsilon_ia_i=\varepsilon_{i+1}a_{i+1}$ and
the common sign $\varepsilon_i$ is $-1$.

The signed major index $\smaj(\a,\varepsilon)$ is the sum
of the signed descents.
\end{definition}

\begin{theorem}
Let 
\begin{equation}
P_n(t,q)=g_n\left(\frac{1-x}{1-q}\right)_{t=-x}=\sum_{i,j}a_n(i,j)t^iq^j\,  
\end{equation}
(where $g_n$ is interpreted as a commutative symmetric function).
Then, $a_n(i,j)$ is equal to the number of signed parking functions
$(\a,\varepsilon)$ of length $n$ with $i$ minus signs  and $j$ signed
inversions, or with signed major index $j$.
In particular, the coefficient of $t^k$ in $P_n(t,0)$ is equal to the number
of nondecreasing signed parking functions  $(\a,\varepsilon)$ with $k$ minus
signs.
\end{theorem}

\Proof The signed inversions (resp. descents) of  $(\a,\varepsilon)$
are the ordinary inversions (resp. descents) of $(\std(\a),\varepsilon)$.
Hence, $\sinv$ and $\smaj$ have the same distribution
over signed parking functions of length $n$.

Now, on the one hand, it has been shown in \cite{NTsuper} that
\begin{equation}\label{qtFsigma}
(q)_n\F_\sigma\left(\frac{\mid 1-x}{1-q\mid }\right)
=\sum_{\varepsilon\in\{\pm 1\}^n}(-x)^{m(\varepsilon)}
                                 q^{\maj(\sigma,\varepsilon)}\,,
\end{equation}
where $m(\varepsilon)$ is the number of minus signs,
and $\maj$ is computed w.r.t. the usual order of $\Z$.

On the other hand, the map $s:\ \PQSym\rightarrow \PQSym^{(2)}$
defined by
\begin{equation}
\F_\a\longmapsto \sum_{\varepsilon\in\{\pm 1\}^n}\F_{(\a,\varepsilon)}
\end{equation}
is a morphism of algebras, and the map
$\phi_{q,x}:\ \PQSym^{(2)}\rightarrow \KK(q)[x]$
defined by
\begin{equation}
\F_{(\a,\varepsilon)} \longmapsto
 (-x)^{m(\varepsilon)}\frac{q^{\smaj(\a,\varepsilon)}}{(q)_n}
\end{equation}
(where $(q)_n=(1-q)\cdots(1-q^n)$)
is a character of $\PQSym^{(2)}$ (also if one replaces $\smaj$
by $\sinv$). 

Thanks to the characterization of signed inversions (or descents) as ordinary
inversions (or descents) of the standardization, we see that
\begin{equation} 
\phi_{q,x}\circ s (\F_\a) = 
\sum_\varepsilon 
 (-x)^{m(\varepsilon)}\frac{q^{\maj(\std(\a),\varepsilon)}}{(q)_n}
= \F_{\std(\a)}\left(\frac{\mid 1-x}{1-q\mid }\right)\,.
\end{equation}
Thus, 
\begin{equation}
\phi_{q,x}\circ s (G)=g\left(\frac{1-x}{1-q}\right)
\end{equation}
where $g$ can now be interpreted as a commutative
symmetric function.  \qed 

For example, in
\begin{equation}
P_2(t,q)=(1+2q)t^2+(3+3q)t+ (2+q) 
\end{equation}
the coefficient of $t^0$ comes from $11,12,21$, the coefficient of $t^1$
from $\bar 11,\bar 12,\bar 21, 1\bar 1 ,1\bar 2,2\bar 1$ and the coefficient
of $t^2$ from $\bar 1\bar 1, \bar 1\bar 2,\bar 2\bar 1$.
Similarly, in 
\begin{equation}
\begin{split}
P_3(t,q)=
& \left( 5\,{q}^{2}+5\,q+5\,{q}^{3}+1 \right) {t}^{3}+
 \left( 10\,{q}^{3}+16\,{q}^{2}+16\,q+6 \right) {t}^{2}\\
&+ \left( 6\,{q}^{3}+16\,{q}^{2}+16\,q+10 \right) t+{q}^{3}+5\,{q}^{2}+5\,q+5,
\end{split}
\end{equation}
the constant term  in $q$ can be explained as follows.
There are 5 parking functions of standardized 123,
3 for 132 and 312, 2 for 213 and 231 and 1 for 321. 
To get an increasing word, 123 can be signed as 123 and $\bar 123$,
contributing 5 times $1+t$. Similarly, 132 does not contribute
and 312,  yields $3$ times $(t+t^2)$ corresponding to $\bar 312$ and $\bar 3 \bar 1 2$.
231 does not contribute and 213 gives twice $(t+t^2)$. Finally, 321
contributes a term $t^2+t^3$ and we end up with
\begin{equation}\begin{split}
5(1+t)+5(t+t^2)+(t^2+t^3)= 5+10t+6t^2+t^3&=(t+1)((t+1)^2+3(t+1)+1)\\
&=(t+1)c_3(t+1)\,.
\end{split}
\end{equation}
The corresponding signed parking functions are given by an encoding
of Schr\"oder paths extending the classical encoding of Dyck paths by
nondecreasing parking functions. A Dyck path is a sequence of steps
$u=(1,1)$ (up) and $d=(1,-1)$ (down), going from $(0,0)$ to $(2n,0)$
without crossing the horizontal axis. Encoding each up step by the
number of its diagonal, we obtain a nondecreasing parking function.
For example, $uuududdudd$ yields $11124$.

\begin{figure}[ht]
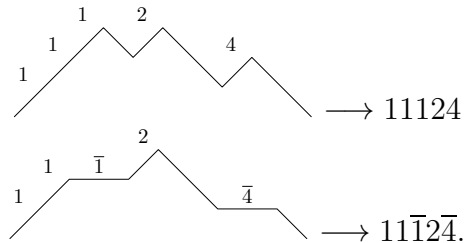

\begin{center}
%  \resizebox{7cm}{!}{\includegraphics{schroeder.pdf}}
\includegraphics[width=4cm]{figs.1}
$\longrightarrow 11124$
\end{center}
\begin{center}
\includegraphics[width=4cm]{figs.2}
$\longrightarrow 11\overline1 2\overline 4$.
\end{center}
\caption{\label{fig:paths}Encoding of Dyck and Schr\"oder paths.}
\end{figure}

Now, Schr\"oder paths are obtained from Dyck paths by replacing
some peaks (factors $ud$) by a horizontal step $h=(2,0)$.
If we replace $k$ by $\bar k$ when a peak with up step in diagonal
$k$ has been replaced by a horizontal step, we obtain an encoding
of Schr\"oder paths by (some) signed nondecreasing parking functions.
Sorting these words according to the order of the signed alphabet,
we obtain precisely those signed parking functions which have no
signed inversion. For example, replacing the first and last peaks
in the previous example by horizontal steps, we obtain the
path $uuhuddhd$, encoded by $11\bar 1 2\bar 4$, whose sorting
is $\bar 4\bar 1112$ (see Figure \ref{fig:paths}).

The coefficients of the polynomials $P_n(t):=P_n(t,0)$ form the sequence
(triangle $T(n,k)$)
\cite[A060693]{Slo}. Since $P_n(-t)$ is the coefficient of $z^n$
in $g((1-t)z)$, unique solution of the functional equation
\begin{equation}
g=\sum_{n\ge 0} S_n(1-t)z^ng^n =1+(1-t)\frac{zg}{1-zg}
\end{equation}
we have the generating series
\begin{equation}
\sum_{n\ge 0}z^n P_n(t)=\frac{1-tz-\sqrt{(1-tz)^2-4z}}{2z}
\end{equation}
which coincides with that of \cite[A060693]{Slo} (up to the
constant term, which is 1 for us).

\begin{corollary}
The coefficient of $t^k$ in $P_n(t)$ is the number of Schr\"oder
paths of semi-length $n$ with $k$ horizontal steps.
In particular, we recover the formula of Lassalle, since it is known that
the generating polynomial for horizontal steps in Schr\"oder paths is
$ (t+1)c_n(t+1)$. Hence
\begin{equation}
tc_n(t)=P_n(t-1)\,.
\end{equation}
\end{corollary}
\qed

%%%%%%%%%%%%%%%%%%%%%%%%%%%%%%%%%%%%%%%%%%%%%%%%%%%%%%%%%%%%%%%%%%%%%%%%%%%%%%%
\subsection{A character of $\SQSym$}

The character yielding the Narayana polynomials has a
simple description at the level of $\SQSym$.
Let $b(\q)$ be the number of bars of a parking quasi-ribbon $\q$.
Then, $\chi(\P_\q)=(1+t)t^{b(\q)}$ is a character of $\SQSym$.

Let again $G(z)=\sum_\q z^{|\q|}\P_\q=\sum_\a z^{|\a|}\F_\a$ be  the sum of
all parking functions (with a homogeneity variable $z$). Then,
\begin{equation}
\chi(G(z)) = g(z(1-x))|_{x=-t}
\end{equation}
so that
\begin{equation}
\chi(G_n) = (1+t)c_n(1+t)\,.
\end{equation}

Thus, $c_n(1+t)$ counts parking quasi-ribbons according to the number of
bars. 
As we have already seen, the classical interpretation of $c_n(1+t)$ is in
terms of Schr\"oder paths (the coefficient of $t^k$ is the number of paths of
semilength $n$ containing exactly $k$ peaks but no peaks at level one
\cite[A126216]{Slo}).
Thus, parking quasi-ribbons with $k$ bars are in bijection with these paths.

A simpler interpretation still in terms of Schr\"oder paths is: the number
of paths with $k$ horizontal steps, and no peak after the last horizontal step.

Replacing the last horizontal step by a peak gives the missing paths, wich
explains the factor $(1+t)$.

%%%%%%%%%%%%%%%%%%%%%%%%%%%%%%%%%%%%%%%%%%%%%%%%%%%%%%%%%%%%%%%%%%%%%%%%%%%%%%%
\subsection{Other characters}

Any character of $\Sym$ (hence also of $Sym$) can be lifted 
to $\PQSym$, thanks to the following property:

\begin{lemma}
The map $\psi: \ \PQSym\rightarrow \Sym$ defined by
\begin{equation}
\psi(\F_\a)=\frac{1}{n!}S^{t(\a)} \quad\text{for $\a\in\PF_n$}
\end{equation}
where $t(\a)$ is the packed evaluation of $\a$,
is a morphism of algebras.
\end{lemma}
\Proof Let $\a\in\PF_n$ and $\b\in\PF_m$, $I=t(\a)$ and $J=t(\b)$.
Then, the product $\F_\a\F_\b$ contains $\frac{(n+m)!}{n!m!}$ terms $\F_\c$,
which all have packed evaluation $t(\c)=I\cdot J$. \qed

Of course, this is not a morphism of coalgebras, as this would
imply by duality an embedding of $QSym$ into a free associative algebra.

As a last illustration, let us consider the classical character
of $Sym$, given by evaluation on a binomial element. Recall that the symmetric
functions of a binomial element $\alpha$  are given by
\begin{equation}
e_n(\alpha)=
\binom{\alpha}{n}\ \text{or, equivalently,}\
h_n(\alpha)=\binom{\alpha+n-1}{n},\ p_n(\alpha)=\alpha\,.
\end{equation}
Thus,
\begin{equation}
\frac1{n!}S^I(\alpha)
=\prod_{k}\frac{1}{i_k!}\cdot\frac{\prod_{k}\prod_{j_k=1}^{i_k}
 (\alpha+j_k-1)}{n!}
=\chi(S^I)\cdot\frac{1}{n!}Z_I(\alpha)
\end{equation}
where $\chi$ is the character of $\Sym$ given by $S_n\mapsto\frac{1}{n!}$, and 
\begin{equation}
Z_I(\alpha) =\sum_{\sigma\in {\rm Stab} (w)}\alpha^{\#{\rm cycles}(\sigma)}
\end{equation}
is the cycle enumerator of the stabilizer of any word  $w$ of (packed)
evaluation $I$.
Dividing by $\chi(S^I)$, we have:

\begin{proposition}
The map $\psi_\alpha:\ \PQSym\rightarrow \K[\alpha]$ defined by
\begin{equation}
\psi_\alpha(\F_\a)=\frac1{n!}Z_{t(\a)}(\alpha)
\end{equation}
is a character of $\PQSym$.
\qed
\end{proposition}

The image of the sum of all parking functions of length $n$ is thus
\begin{equation}
\psi_\alpha(G_n)=g_n(\alpha)=\frac1{n+1}h_n((n+1)\alpha)
=\frac1{n+1} \binom{(n+1)\alpha+n-1}{n}
=\frac{P_n(\alpha)}{n!}
\end{equation}
where 
\begin{equation}
P_n(\alpha)=\alpha\prod_{k=1}^{n-1}((n+1)\alpha+k)\,.
\end{equation}

From the previous considerations, we obtain a combinatorial
interpretation of the coefficients of $P_n(\alpha)$:

\begin{corollary}
The coefficient of $\alpha^k$ in $P_n(\alpha)$ is equal to
the number of pairs $(\a,\sigma)$ where $\a\in\PF_n$ and
$\sigma$ is a permutation in $\SG_n$
with $k$ cycles,  such that $\a\sigma=\a$.
\qed
\end{corollary}

The first polynomials $P_n$ are

\begin{eqnarray}
P_1(\alpha)&=&\alpha,\\
P_2(\alpha)&=&3\alpha^2+\alpha,\\
P_3(\alpha)&=&16\alpha^3+12\alpha^2+2\alpha,\\
P_4(\alpha)&=&125\alpha^4+150\alpha^3+55\alpha^2+6\alpha.
\end{eqnarray}

We can observe a curious property of the reciprocal
polynomials evaluated at $q-1$: if we set
\begin{equation}
Q_n(q)=(q-1)^{n-1}P_n\left(\frac1{q-1} \right) =\prod_{k=2}^n((n+1-k)q+k)
\end{equation}
we obtain a $q$-analogue of $(n+1)^{n-1}$ whose coefficients form the triangle

\begin{equation}
\begin{matrix}
1 &&&\\
2 &1&&\\
6& 8& 2&\\
24 & 58& 37& 6
\end{matrix}
\end{equation}
The first column is clearly $n!$. The second column is sequence A002538 in the
OEIS \cite{Slo}, which counts the number of edges in the Hasse diagram of the
Bruhat order on $\SG_n$. Since the first column is its number of vertices, one
may wonder whether the other numbers also have an interpretation in this
context.
This will be left as an open problem.

%%%%%%%%%%%%%%%%%%%%%%%%%%%%%%%%%%%%%%%%%%%%%%%%%%%%%%%%%%%%%%%%%%%%%%%%%%%%%%%
%%%%%%%%%%%%%%%%%%%%%%%%%%%%%%%%%%%%%%%%%%%%%%%%%%%%%%%%%%%%%%%%%%%%%%%%%%%%%%%
%%%%%%%%%%%%%%%%%%%%%%%%%%%%%%%%%%%%%%%%%%%%%%%%%%%%%%%%%%%%%%%%%%%%%%%%%%%%%%%
%%%%%%%%%%%%%%%%%%%%%%%%%%%%%%%%%%%%%%%%%%%%%%%%%%%%%%%%%%%%%%%%%%%%%%%%%%%%%%%

\end{document}